\documentclass[a4paper,12pt,leqno]{article}
\usepackage{latexsym}
\usepackage[all]{xy}

\usepackage{amssymb} 
\usepackage{amsmath} 
\usepackage{theorem}

\def\Z{{\mathbb{Z}}}
\def\K{{\mathbb{K}}}
\def\R{{\mathbb{R}}}
\def\C{{\mathbb{C}}}
\def\A{{\mathcal{A}}}
\def\B{{\mathcal{B}}}
\def\C{{\mathcal{C}}}

\DeclareMathOperator{\codim}{codim}

\DeclareMathOperator{\Der}{Der}

\DeclareMathOperator{\Shi}{Shi}

\numberwithin{equation}{section}

\newcommand{\owari}{\hfill$\square$}

\theoremstyle{break}
\newtheorem{theorem}{Theorem}[section]
\newtheorem{prop}[theorem]{Proposition}
\newtheorem{cor}[theorem]{Corollary}

\newtheorem{define}[theorem]{Definition}

\newtheorem{rem}[theorem]{Remark}

\newtheorem{example}[theorem]{Example}
\newtheorem{problem}[theorem]{Problem}
\newtheorem{conj}[theorem]{Conjecture}

\title{Restrictions of free arrangements and 
the division theorem}
\author{
Takuro Abe\footnote
{
Institute of Mathematics for Industry, 
Kyushu University,
Fukuoka 819-0395, Japan.
email:abe@imi.kyushu-u.ac.jp
\ \ 
\textit{MSC primary}. 32S22, 52S35.
}
}
\date{\today}

\begin{document}

\maketitle

\begin{abstract}
This is a survey and research note on the modified Orlik conjecture derived from 
the division theorem introduced in 
\cite{A2}. The division theorem is 
a generalization of classical addition-deletion theorems for free arrangements. The division 
theorem can be regarded as a modified converse of the Orlik's conjecture 
with a combinatorial condition, i.e., 
an arrangement is free if the restriction is free and the characteristic polynomial of 
the restriction divides that of an arrangement. 
In this article we recall, summarize, pose and re-formulate some of results and 
problems related to the division theorem based on \cite{A2}, and study the modified Orlik's conjecture with partial answers.
\end{abstract}

\section{Introduction}
Let $\K$ be an arbitrary field, $V=\K^\ell$ and 
$S=\K[x_1,\ldots,x_\ell]$ the coordinate ring of $V^
*$. Let $\A$ be an arrangement of hyperplanes in $V$, i.e.,
a finite collection of linear hyperplanes in $V$. For $H \in \A$ fix 
a linear form $\alpha_H \in V^*$ such that $\mbox{ker}(\alpha_H)=
H$. For $\Der S:=\oplus_{i=1}^\ell S \partial_{x_i}$, a logarithmic derivation module 
$D(\A)$ of $\A$ is defined by 
$$
D(\A):=\{\theta \in \Der S \mid 
\theta(\alpha_H) \in S \alpha_H\ (\forall H \in \A)\}.
$$
$D(\A)$ is a reflexive $S$-module, and not free in general. We say $\A$ is free with
exponents $\exp(\A)=(d_1,\ldots,d_\ell)$ if there is derivations $\theta_1,\ldots,\theta_\ell 
\in D(\A)$ such that $D(\A)=
\oplus_{i=1}^\ell S \theta_i$ and $\deg \theta_i(\beta_i)=d_i\ (i=1,\ldots,\ell)$ for some 
linear form $\beta_1,\ldots,\beta_\ell$ such that $\theta_i(\beta_i) \neq 0$.

Let $L(\A):=\{\cap_{H \in \B} H\mid \B \subset \A\}$ be the intersection lattice of $\A$, and 
$L_i(\A):=\{X \in L(\A) \mid \codim_VX=i\}$. A flag of $\A$ is a set 
$\{X_i\}_{i=0}^{\ell}$ such that $X_0 \subset \cdots \subset X_\ell$ and 
$X_i \in L_i(\A)$. The M\"{o}bius function $\mu:L(\A) \rightarrow 
\Z$ is defined by, $\mu(V
)=1$ and by $\mu(X):=-\sum_{X \subsetneq Y \subset V} \mu(Y)$ for $X \neq V$. The Poincar\`{e} polynomial $\pi(\A;t)$ of $\A$ is defined by 
$\pi(\A;t):=\sum_{X \in L(\A)} \mu(X)(-t)^{\codim_V X}$. Also, the 
characteristic polynomial $\chi(\A;t)$ of $\A$ is 
defined by $\chi(\A;t):=t^\ell \pi(\A;-t^{-1})$. 
It is known that 
$\pi(\A;t)$ coincides with the topological Poincar\`{e} polynomial of the complement 
$M(\A):=V \setminus \cup_{H\in \A} H$ when $\K=\mathbb{C}$. Hence 
the coefficient $b_i(\A)$ of $t^i$ in $\pi(\A;t)$ is nothing but the 
$i$-th Betti number of $M(\A)$. For $X \in L(\A)$, a localization $\A_X$ of $\A$ at $X$ is defined by 
$\A_X:=\{H \in \A \mid X \subset H\}$, and the restriction 
$\A^X$ of $\A$ onto $X$ is defined by $\A^X:=\{H \cap X\mid
H \in \A \setminus \A_X\}$. Note that $\A_X$ is an arrangement in $V$, but 
$\A^X$ is that in $X\simeq \K^{\dim X}$.

Free arrangements have been intensively studied by several mathematicians, 
and that research has been the most important among the study of algebraic aspects of 
an arrangement. 
To check the freeness of given arrangement, or 
to construct a new free arrangement is very difficult though that is very fundamental. 
For that purpose, Terao's addition-deletion and 
restriction theorems have been the most useful and important.

\begin{theorem}[\cite{T1}, Addition-deletion and restriction theorems]
For $H \in \A$, let $\A':=\A \setminus \{H\}$. Then for the triple $(\A,\A',\A^H)$, two 
of the following three implies the third:
\begin{itemize}
\item[(1)]
$\A$ is free with $\exp(\A)=(d_1,\ldots,d_{\ell-1},d_\ell)$.
\item[(2)]
$\A'$ is free with $\exp(\A')=(d_1,\ldots,d_{\ell-1},d_\ell-1)$.
\item[(3)]
$\A^H$ is free with $\exp(\A^H)=(d_1,\ldots,d_{\ell-1})$.
\end{itemize}
Moreover, all the three above hold if both $\A$ and 
$\A'$ are free.
\label{adddel}
\end{theorem}

In \cite{A2} the division theorem for free arrangements 
was introduced, which is a generalization of Terao's addition-deletion theorem \ref{adddel}.

\begin{theorem}[Theorem 1.1, \cite{A2}, Division theorem]
$\A$ is free if $\A^H$ is free and $\chi(\A^H;t)$ divides 
$\chi(\A;t)$ for some $H \in \A$.
\label{division}
\end{theorem}

Theorem \ref{division} can be regarded as a converse of modified Orlik's 
conjecture. Orlik's conjecture asserted that 
$\A^H$ is free if $\A$ is free, the counter example to which was 
found by Edelman and Reiner in \cite{ER}. Theorem \ref{division} is a converse of this conjecture 
with one more condition that $\chi(\A^H;t)$ divides $\chi(\A;t)$. Then it is natural to ask 
whether this modified Orlik's conjecture is true or not.

\begin{problem}[Modified Orlik's Conjecture]
Let $\A$ be an $\ell$-arragement and $H \in \A$. Assume that 
$\A$ is free and $\pi(\A^H;t)$ divides $\pi(\A;t)$. (Equivalently, 
$\pi(\A\setminus \{H
\};t)$ divides $\pi(\A;t)$.) Then is $\A^H$ (and hence $\A \setminus \{H\}$) a free 
arrangement?
\label{Orlik}
\end{problem}

It seems that what is stated in Problem \ref{Orlik} is too strong, hence 
we believe that there will be a counter example to Problem \ref{Orlik}. In other words, 
we believe the following conjecture.

\begin{conj}
There exists an arrangement $\A$ and $H \in \A$ such that, for the 
triple $(\A,\A',\A^H)$, it holds that 
\begin{itemize}
\item[(1)]
$\pi(\A^H;t)=\prod_{i=1}^{\ell-1} (1+d_it)$ with $d_1,\ldots,d_{\ell-1} \in \Z$,
\item[(2)]
$\pi(\A^H;t)$ divides
both $\pi(\A;t)$ and $\pi(\A';t)$, and 
\item[(3)]
Neither $\A$ nor $\A'$ are free (or, $\A$ is free and $\A'$ is not free).
\end{itemize}
\label{conj1}
\end{conj}

By Theorem \ref{division}, there are no triples as in Conjecture \ref{conj1} if 
the condition (3) in Conjecture \ref{conj1} holds. Also, there are no such triple 
when $\ell \le 3$ due to \cite{A}. Hence to show Conjecture \ref{conj1}, the assumption 
$\ell \ge 4$ is essential.
However, and surprisingly, 
there have been no example as in Conjecture \ref{conj1}.



The purpose of this article is to consider in which condition Problem \ref{Orlik} is true. 
The key role is played by the second Betti number $b_2(\A)$ of the complement 
of $\A$ when $\K=\mathbb{C}$. Namely, $b_2(\A)$ is the coefficient of 
$t^{\ell-2}$ of $\chi(\A;t)$. 
One of the answer is the following, which is a main result in this article.

\begin{theorem}
Assume that $\A$ is free and $b_2(\A)=b_2(\A^H)+(|\A|-|\A^H|)|\A^H|$. 
Then $\A^H$ is free if 
there is $L \in \A \setminus \{H\}$ such that $\A \setminus \{L\}$ is free and 
$|\A_{L \cap H}| \ge 3$.
\label{deform}
\end{theorem}

Note that the equation 
$b_2(\A)=b_2(\A^H)+(|\A|-|\A^H|)|\A^H|$ holds when $\chi(\A^H;t)$ divides 
$\chi(\A;t)$. What is interesting in Theorem \ref{deform} is, to determine the freeness of 
the restriction, a freeness of some other restriction works.

The other main result in this article is to give an easy sufficient condition for 
$\A'$ not to be free even when $\A$ is free. This gives us an easy sufficient 
condition for the equation on the second Betti numbers above
not to be true.

\begin{theorem}
Assume that $\A$ is a free $\ell$-arrangement with $\ell \ge 3$. 
Let $H \in \A$ and $\A':=\A \setminus \{H\}$. Then $\A'$ is not free if 
there is $X \in L_2(\A^H)$ such that one of the following three holds:
\begin{itemize}
\item[(1)]
$\A_X^H=\{K_1,K_2\}$ and $m^H(K_1) > 1,\ m^H(K_2) > 1$.
\item[(2)]
$\A_X^H=\{K_1,K_2,K_3\}$ and $m^H(K_i) \ge 2$ for $i=1,2,3$.
\item[(3)]
$\A_X^H=\{K_1,K_2,K_3\}$ and $m^H(K_1) \ge 3, \ 
m^H(K_2) \ge 2$.
\end{itemize}
Here $m^H:\A^H \rightarrow \Z_{>0}$ is the Ziegler multiplicity on $\A^H$ defined by 
$$
m^H(K):=|\{L  \in \A' \mid H \cap L=K\}|
$$
for $K \in \A^H$.
In particular, if one of the three above holds, then it holds that 
$b_2(\A)>b_2(\A^H)+(|\A|-|\A^H|)|\A^H|$.
\label{notfree}
\end{theorem}

Hence the arrangement $\A$ in Problem \ref{Orlik} has a very special geometry. 
Let us show an application of Theorem \ref{notfree}.

\begin{example}
Let $\A$ be an arrangement in $V=\R^4$ defined by 
$$
\prod_{i=1}^4 x_i \prod_{1 \le i < j \le 4,\ (i,j) \neq (3,4)} (x_i-x_j)=0.
$$
Then it is easy to see that $\A$ is free with $\exp(\A)=(1,2,3,3)$. Let 
$H:=\{x_1=0\} \in \A$, and show that $\A \setminus \{H\}$ is not free. 
The Ziegler restriction $(\A^H,m^H)$ of $\A$ onto $H$ is defined by 
$$
(\prod_{i=2}^4 x_i^2 )(x_2-x_3)(x_2-x_4)=0.
$$
Let $X:=\{x_3=x_4=0\} \in L_2(\A^H)$. Then $\A_X^H=\{
x_3^2x_4^2=0\}$, which satisfies the condition (1) in Theorem \ref{notfree}. Hence 
$b_2(\A) > b_2(\A^H)+(|\A|-|\A^H|)|\A^H|$ and 
$\A \setminus 
\{H\}$ is not free.
\label{A4}
\end{example}




The organization of this article is as follows. In \S 2 we recall several results and definitions for the 
proof. This section contains some re-formulation of results in \cite{A2}. 
In \S 3, first we give some partial answers to Problem \ref{Orlik} which follows immediately from 
the division theorem and other results in \cite{A2}. After that, we show 
Theorems \ref{deform} and \ref{notfree}. In \S 4 we observe the similarity of supersolvable and 
divisionally free arrangements.
\medskip

\noindent
\textbf{Acknowledgments}. 
The author is grateful to the referee for the careful reading of this 
paper with a lot of important comments. 
This work is partially supported by 
JSPS Grants-in-Aid for Young Scientists
(B)
No. 24740012.

\section{Preliminaries}

In this section let us recall several results we will use for the proof of main results. 
The first one is the most important result among the theory of free arrangements.

\begin{theorem}[\cite{T2}, Terao's factorization]
Assume that $\A$ is free with $\exp(\A)=(d_1.\ldots,d_\ell)$. Then 
$\chi(\A;t)=\prod_{i=1}^\ell (t-d_i)$. In particular, 
$\A$ is not free if $\chi(\A;t)$ is irreducible over $\Z$.
\label{factorization}
\end{theorem}

Next let us recall some fundamental definitions and results for multiarrangements. 
For an arrangement $\A$, let $m:\A \rightarrow \Z_{\ge 0}$ be a multiplicity function. Then the pair 
$(\A,m)$ is called a multiarrangement, and we can define the logarithmic derivation module 
$D(\A,m)$ by
$$
D(\A,m):=\{
\theta \in \Der S \mid \theta(\alpha_H) \in S \alpha_H^{m(H)}\ (\forall H \in \A)\}.
$$
Then we can define the freeness and exponents for multiarrangements in the same manner 
as for arrangements. 

From an arrangement, we may define a multiarrangement canonically. 
For an arrangement $\A$ and $H \in \A$, define a multiarrangement $(\A^H,m^H)$, called the 
Ziegler restriction of $\A$ onto $H$, by 
$m^H(X):=|\A_X|-1$ for $X \in \A^H$. Then the following is the most fundamental.

\begin{theorem}[\cite{Z}]
Assume that $\A$ is free with $\exp(\A)=(1,d_2,\ldots,d_\ell)$. Then $(\A^H,m^H)$ is free with 
$\exp(\A^H,m^H)=(d_2,\ldots,d_\ell)$.
\label{Ziegler}
\end{theorem}

The next result is a generalization of Theorem \ref{adddel} to multiarrangements. For 
details and definitions on the Euler restriction, see \cite{ATW2}.

\begin{theorem}[\cite{ATW2}, Theorem 0.8]
Let $(\A,m)$ be a multiarrangment, $H \in \A$ and let 
$\delta_H:\A \rightarrow \{0,1\}
$ be a multiplicity such that $\delta_H(L)=1$ only when $H=L$. Then any two of the 
following imply the third:
\begin{itemize}
\item[(1)]
$(\A,m)$ is free with $\exp(\A,m)=(d_1,\ldots,d_{\ell-1},d_\ell)$.
\item[(2)]
$(\A,m-\delta_H)$ is free with $\exp(\A,m-\delta_H)=(d_1,\ldots,d_{\ell-1},d_\ell-1)$.
\item[(3)]
$(\A^H,m^*)$ is free with $\exp(\A^H,m^*)=(d_1,\ldots,d_{\ell-1})$,
\end{itemize}
where $(\A^H,m^*)$ is the Euler restriction of $(\A,m)$ onto $H \in \A$. 
Moreover, all the three above hold if both $(\A,m)$ and 
$(\A,m-\delta_H)$ are free.
\label{multiadddel}
\end{theorem}

The following is a freeness criterion by using the second Betti number and 
the Ziegler restriction. For details, see \cite{AY}. Also, for the definiton of the 
second Betti number of a multiarrangement, see \cite{ATW}.

\begin{theorem}[\cite{AY}, Theorem 5.1]
Let $\A$ be a central $\ell$-arrangement, $H \in \A$ and 
$(\A^H,m^H)$ the Ziegler restriction of $\A$ onto $H$. Then 
$\A$ is free if and only if $(\A^H,m^H)$ is free and 
$b_2(\A)=|\A|-1+b_2(\A^H,m^H)$. In particular, 
$b_2(\A) \ge |\A|-1+b_2(\A^H,m^H)$.
\label{lfree}
\end{theorem}



Let us introduce two more results
from \cite{A2}. Since the formulations of these results are slight different from those 
in the original version in \cite{A2}, we give proofs for the completeness.

The first one
is the following proposition, which says that 
the Ziegler and Euler restriction commutes if 
there is a division $\chi(\A^H;t) \mid \chi(\A;t)$.

\begin{prop}[\cite{A2}, cf. Theorem 1.7]
Let $\A$ be an $\ell$-arrangement, $H \in \A$ and 
$(\A^H,m^H)$ be the Ziegler restriction of $\A$ onto $H$. 
Let $X \in \A^H$ with $m^H(X) \ge 2$.
Assume that 
$b_2(\A)=b_2(\A^H)+(|\A|-|\A^H|)|\A^H|$. Then 
\begin{itemize}
\item[(1)]
the Ziegler restriction of $\A^H$ onto $X$
coincides with 
the Euler restriction of $(\A^H,m^H)$ onto $X$, and 
\item[(2)]
$b_2(\A^H,m^H)-b_2(\A^H,m^H-\delta_X)=|\A^H|-1$.
\end{itemize}
\label{key1}
\end{prop}

\noindent
\textbf{Proof}. Immediate from 
Theorem 1.7 and its proof in \cite{A2}.\owari
\medskip


\begin{prop}[$(b_1,b_2)$-inequality, cf., \cite{A2}, Corollary 4.10 ]
Let $H \in \A$. Then 
$$
b_2(\A) \ge b_2(\A^H)+(|\A|-|\A^H|)|\A^H|.
$$
In particular, 
$$
b_2(\A) \ge \sum_{i=0}^{\ell-2}
(|\A^{X_i}|-|
\A^{X_{i+1}}|)|\A^{X_{i+1}}|
$$
for any flag $\{X_i\}_{i=0}^{\ell-1}$ of $\A$.
\label{b1b2}
\end{prop}

\noindent
\textbf{Proof}.
Let $b_2(d\A)$ denote the coefficient of $t^2$ in 
$\pi(\A;t)/(1+t)$. Then the equation (4.1) in \cite{A2} shows that 
$$
b_2(d\A) 
\ge b_2(d\A^H)+(|\A|-|\A^H|)(|\A^H|-1).
$$
Since $b_2(d\A)+|\A|-1=b_2(\A)$ and 
$|\A|-1=(|\A|-|\A^H|)+(|\A^H|-1)$, we have 
$b_2(\A) \ge b_2(\A^H)+
(|\A|-|\A^H|)|\A^H|$. The next inequality holds by applying the same argument to the 
Corollary 4.10 in \cite{A2}. \owari
\medskip

Here we give an observation. As in Theorem \ref{adddel}, the freeness of 
$\A$ and $\A'$ implies the freeness of each member of the triple. 
Here we do not have to consider 
the freeness of $\A^H$. Then, when each member of the triple 
is free if we assume the freeness of $\A^H$? The answer is immediate from 
Theorem \ref{division}.

\begin{prop}
Let $\A$ be an arrangement and $(\A,\A',\A^H)$ the triple 
with respect to $H \in \A$. Then 
each member of the triple 
is free if and only if $\A^H$ is free and 
$b_2(\A)=b_2(\A^H)+(|\A|-|\A^H|)|\A^H|$.
\label{main2}
\end{prop}

\noindent
\textbf{Proof}. The ``only if'' part is nothing but Theorem \ref{adddel}. 
The ``if'' part follows also immediately by Theorem \ref{division}.\owari
\medskip

\section{A partial results and the proof of main results}

Before the proof of Theorem \ref{deform}, 
let us give some partial answer which follows immediately from 
the division theorem and $(b_1,b_2)$-inequality.


\begin{theorem}
Let $\A$ be a free arrangement with $\exp(\A)=(d_1,\ldots,d_\ell)$. 
Take $H \in \A$ and $X \in \A^H$ such that 
\begin{itemize}
\item[(1)] 
$\A^X$ is free with $\exp(\A^X)=(d_1,\ldots,d_{\ell-2})$, and 
\item[(2)]
$|\A|-|\A^H|=d_\ell,\ (\mbox{hence automatically, }|\A^H|-|\A^X|=d_{\ell-1}$).
\end{itemize}
Then $\A^H$ is also free with $\exp(\A^H)=(d_1,\ldots,d_{\ell-1})$.
\label{sandwich}
\end{theorem}


\noindent
\textbf{Proof}.
By Theorem \ref{division}, it suffices to show that $b_2(\A^H)=
\sum_{1 \le i < j \le \ell-1} d_i d_j$. By Proposition \ref{b1b2}, it holds that 
\begin{eqnarray*}
b_2(\A)&=&\sum_{1 \le i < j \le \ell} d_i d_j \ge b_2(\A^H) +d_\ell (d_1+\cdots+d_{\ell-1}),\ \mbox{and}\\
b_2(\A^H)  &\ge& b_2(\A^X) +d_{\ell-1} (d_1+\cdots+d_{\ell-2})=
\sum_{1 \le i < j \le \ell-2} d_i d_j+d_{\ell-1} (d_1+\cdots+d_{\ell-2}).
\end{eqnarray*}
Hence 
$$
\sum_{1 \le i < j \le \ell-1} d_i d_j \le b_2(\A^H) \le \sum_{1 \le i < j \le \ell-1}d_i d_j.
$$
Hence 
$b_2(\A^H)=b_2(\A^X)+(|\A^H|-|\A^X|)|\A^X|$, and Theorem \ref{division} completes the proof. \owari
\medskip

In Theorem \ref{sandwich}, 
we apply the proof of Theorem \ref{division} conversely. Hence Theorem \ref{sandwich} may be 
regarded as an application of the proof of Theorem \ref{division} and 
Proposition \ref{b1b2}. 
A useful part of Theorem \ref{sandwich} is, if we know the exponents of $\A$ and $\A^{
X}$, then we can check the freeness between them just by computing the number of hyperplanes 
in it (we do not need any information on the second Betti number!). Hence practically, or when we want to 
check some hereditary freeness (see \cite{OT2}), Theorem \ref{sandwich} and the following 
corollaries are useful.


\begin{cor}
Let $\A$ be a free arrangement with $\exp(\A)=(d_1,\ldots,d_\ell)$. 
Take $X_i \in L_i(\A)\ (i=1,\ldots,k)$ with $X_1 \supset \cdots \supset X_k$ such  that 
\begin{itemize}
\item[(1)] 
$\A^{X_k}$ is free with $\exp(\A^{X_k})=(d_1,\ldots,d_{\ell-k})$.
\item[(2)]
$|\A^{X_i}|-|\A^{X_{i+1}}|=d_{\ell-i}\ (i=0,\ldots,k-1)$.
\end{itemize}
Then $\A^{X_i}$ is also free with $\exp(\A^{X_i})=(d_1,\ldots,d_{\ell-i})$ for 
$i=1,\ldots,k-1$.
\label{mille}
\end{cor}

\noindent
\textbf{Proof}.
By Proposition \ref{b1b2}, it holds that 
$$b_2(\A^{X_{k-1}}) \ge 
(|\A^{X_{k-1}}|-|\A^{X_{k}}|
)|\A^{X_{k}}|+b_2(\A^{X_{k}})=
\sum_{1 \le i<j \le \ell-k+1} d_id_j.
$$
On the other hand, again by applying $(b_1,b_2)$-inequality, we have 
\begin{eqnarray*}
\sum_{1 \le i <j\le \ell}d_jd_i&=&
b_2(\A) \\
&\ge&
b_2(\A^{X_1})+(|\A|-|\A^{X_{1}}|
)|\A^{X_{1}}|=b_2(\A^{X_1})+d_\ell\sum_{i=1}^{\ell-1}d_i\\
&\ge&
b_2(\A^{X_2})+(|\A^{X_1}|-|\A^{X_{2}}|
)|\A^{X_{2}}|+d_\ell\sum_{i=1}^{\ell-1}d_i\\
&\ge&\cdots\\
&\ge&
b_2(\A^{X_{k-1}})+(\sum_{i=\ell-k+2}^\ell d_i) (
\sum_{j=1}^{i-1}d_j).
\end{eqnarray*}
Hence it holds that 
$$
b_2(\A^{X_{k-1}}) \le 
\sum_{1 \le i  <j\le \ell}d_jd_i-(\sum_{i=\ell-k+2}^\ell d_i) (
\sum_{j=1}^{i-1}d_j)=
\sum_{1 \le i  <j\le \ell-k+1}d_id_j.
$$
Combine these two inequalities to obtain 
$$
b_2(\A^{X_{k-1}})=\sum_{1 \le i <j \le \ell-k+1} d_i d_j.
$$
Hence Theorem \ref{division} shows that $\A^{X_{k-1}}$ is free. 
Apply the same argument to all $\A^{X_1},\ldots,\A^{X_{k-2}}$ to complete the proof.
\owari
\medskip

Moreover, we do not need to assume the freeness of $\A$ as follows:

\begin{cor}
Let $\A$ be an arrangement with 
$|\A|=b_1(\A)=d_1+\cdots+d_\ell,\ b_2(\A)=\sum_{1 \le i < j \le \ell} d_i d_j$ for some 
positive integers $d_1,\ldots,d_\ell$. 
Take $X_i \in L_i(\A)\ (i=1,\ldots,k)$ with $X_1 \supset \cdots \supset X_k$ such  that 
\begin{itemize}
\item[(1)] 
$\A^{X_k}$ is free with $\exp(\A^{X_k})=(d_1,\ldots,d_{\ell-k})$.
\item[(2)]
$|\A^{X_i}|-|\A^{X_{i+1}}|=d_{\ell-i}\ (i=0,\ldots,k-1)$.
\end{itemize}
Then $\A^{X_i}$ is also free with $\exp(\A^{X_i})=(d_1,\ldots,d_{\ell-i})$ for 
$i=0,\ldots,k-1$. In particular, we do not need the freeness of $\A^{X_k}$ if 
$k=\ell-2$.
\label{mille2}
\end{cor}

\noindent
\textbf{Proof}. 
Apply the same argument as in the proof of Corollary \ref{mille} repeatedly.
When $k=\ell-2$, this is nothing but the divisonal freeness in Definition 1.5 in \cite{A2} 
(see 
also Definition \ref{DF}).
\medskip

Now let us prove Theorems \ref{deform}.
\medskip

\noindent
\textbf{Proof of Theorem \ref{deform}}.
Let $X:=L \cap H$. Let $\exp(\A)=(d_1,\ldots,d_\ell)$ with $d_1
=1$. Also, let $\exp(\A')=(d_1,\ldots,d_{\ell-1},d_\ell-1)$, where 
$\A':=\A \setminus \{L\}$. Then $\A^L$ is 
free with $\exp(\A^L)=(d_1,\ldots,d_{\ell-1})$ by Theorem \ref{adddel}. 
Hence Theorem \ref{Ziegler} shows that 
both $(\A^H,m^H)$ and $(\A^H,m^H-\delta_X)$ are free with exponents 
$(d_2,\ldots,d_\ell)$ and $(d_2,\ldots,d_{\ell}-1)$ respectively, where 
$\delta_X:\A^H \rightarrow \{0,1\}$ is a multiplicity such that 
$\delta_X^{-1}(1)=X \in \A^H$. Then 
Theorem \ref{multiadddel} shows that the Euler restriction $(\A^X, m^*)$ of 
$(\A^H,m^H)$ onto $X$ is also free with $\exp(\A^X,m^*)=(d_2,\ldots,d_{\ell-1})$ by Theorem \ref{Ziegler}.  
Now recall that $m^H(X) \ge 2$ by the fact that $|\A_X|\ge 3$. Hence Proposition 
\ref{key1} and the equality $b_2(\A)=b_2(\A^H)+(|\A|-|\A^H|)|\A^H|$ show that 
$(\A^X,m^*)=(\A^X,m^X)$,  where $(\A^X,m^X)$ is the Ziegler restriction of $\A^H$ onto $X$. 
Hence $(\A^X,m^X)$ is also free with $\exp(\A^X,m^X)=(d_2,\ldots,d_{\ell-1})$. In 
particular, $|m^X|=d_2+\cdots+d_{\ell-1}=|\A|-d_\ell-1$. On the other hand, again 
Proposition \ref{key1} shows that $|m^*|=
b_2(\A^H,m^H)-b_2(\A^H,m^H-\delta_X)=|\A^
H|-1$. Since $|m^*|=|m^X|$, we have 
$|\A|-d_\ell-1=|\A^H|-1$. Hence $|\A|-|\A^H|=d_\ell$, and 
$|\A^H|=\sum_{i=1}^\ell d_i-d_\ell=d_1+\cdots+d_{\ell-1}$. So 
the equation $b_2(\A)=b_2(\A^H)+(|\A|-|\A^H|)|\A^H|$ shows that 
$b_2(\A^H)=\sum_{1 \le i <j \le \ell-1}d_i d_j$. In particular, 
the coefficient of $t^2$ in $\pi_0(\A^H;t)$ is 
$\sum_{2 \le i< j \le \ell-1} d_i d_j=b_2(\A^X,m^X)$. Hence Theorem \ref{lfree} 
shows that $\A^H$ is free with $\exp(\A^H)=(d_1,\ldots,d_{\ell-1})$. \owari
\medskip

Theorem \ref{deform} has the following corollary.

\begin{cor}
Let $\A$ be an $\ell$-arrangement and 
$H_1,\ldots,H_s \in \A$ be distinct hyperplanes such that 
$\codim \cap_{i=1}^s H_i=2$. 
Assume that 
$b_2(\A)=b_2(\A^{
H_i})+(|\A|-|\A^{H_i}|)|\A^{H_i}|$ for
$i=2,\ldots,s$. Then all the $\A'_i:=
\A \setminus \{H_i\}\ (i=2,\ldots,s)$ and $\B:=\A \setminus \{H_1,\ldots,H_s\}$ 
are free if $H_1$ satisfies the following 
conditions:
\begin{itemize}
\item[(1)]
$\A$ and $\A \setminus \{H_1\}$ are free, and 
\item[(2)]
$|\A_{X}| \ge s+1$ for $X:=\cap_{i=1}^s H_i$.
\end{itemize}
\label{descent2}
\end{cor}

\noindent
\textbf{Proof}. Apply Theorem \ref{deform} to 
each pair $H_1,H_i$ to obtain the statement in Theorem \ref{deform}. 
Hence $\A^{H_i}$ is free for each $i \ge 2$ with exponents$(d_1,d_2,\ldots,d_{\ell-1})$, 
here we assume that $\exp(\A)=(d_1.\ldots,d_\ell),\ 
d_1=1$ and $|\A|-|\A^{H_1}|=d_\ell$. 
Since $|\A_X|\ge s+1$, it holds that 
$$
(\A \setminus \{H_1,\ldots,H_i\})^{H_{i+1}}=\A^{H_i}.
$$
Hence Theorem \ref{adddel} shows that $\A \setminus \{H_1,\ldots,H_i\}$ is free with 
exponents $(1,d_2,\ldots,d_{\ell-1},d_\ell-i)$, which completes the proof, \owari
\medskip

As we can see from the proof of Theorem \ref{deform}, 
the following general fact holds.

\begin{cor}
Assume that $\A$ is free and $\chi(\A^H;t)$ divides $\chi(\A;t)$ for $H \in \A$. Then 
$\A^H$ is free if 
there is $X \in \A^H$ such that $m^H(X) \ge 2$ and
$(\A^H,m^H-\delta_X)$ is free. 
\label{multideform}
\end{cor}

\noindent
\textbf{Proof}. Immediate from 
the proof of Theorem \ref{deform}. \owari
\medskip

\noindent
\textbf{Proof of Theorem \ref{notfree}}.
Since the proof for each given condition is the same, we show only the case (1).
Assume that $\A'$ is free. Then Theorem \ref{adddel} shows that 
$\A^H$ is free with $\exp(\A^H) \subset \exp(\A)$. In particular, the equation 
$b_2(d\A)=b_2(\A^H)+(|\A|-|\A^H|-1)(|\A^H|-1)$ holds. By the proof of Theorem \ref{division}
in \cite{A2}, for every multiplicity $m:\A \rightarrow \Z_{>0}$ with 
$m(Y) \le m^H(Y)\ (Y \in \A^H)$, it holds that 
$b_2(\A^H,m)=b_2(\A^H)+(|m|-|\A^H|)(|\A^H|-1)$. In particular, for any $Y \in \A^H$ with 
$m(Y) \ge 2$, it holds that $b_2(\A^H,m)-b_2(\A^H,m-\delta_Y)=|\A^H|-1$. 

By assumption (1), we may pick a multiplicity $m$ such that $m(K_1)=m(K_2)=2$. 
Let $m^*$ be the Euler multiplicity of $(\A^H,m)$ onto $K_1$. Then it follows that 
$m^*(X)=2$. Hence 
\begin{eqnarray*}
b_2(\A^H,m)-b_2(\A^H,m-\delta_{K_1})&=&
\sum_{Y \in (\A^H)^{K_1}} (b_2(\A^H_Y,m)-b_2(\A^H_Y,m-\delta_{K_1}))\\
&=&\sum_{Y \in (\A^H)^{K_1}} m^*(Y)\\
&>& \sum_{Y \in (\A^H)^{K_1}} (|\A_Y^H|-1)\\
&=& |\A^H|-1
\end{eqnarray*}
by Lemma 3.3 (2), Lemma 3.4  and the assumption
that $m^*(X) =2>1$, which is a contradiction. 
For other cases, use the same argument with the result in \cite{W}.\owari
\medskip

Let us 
see an example how to apply Corollary \ref{mille}.

\begin{example}
Let $\A$ be an arrangement in $V=\R^6$ defined by 
$$
\prod_{1 \le i < j \le 6}(x_i^2-x_j^2)=0.
$$
This is the Weyl arrangement of the type $D_6$, hence free with $\exp(\A)=(1,3,5,5,7,9)$. 
In general, to investigate the freeness of restrictions is very difficult. In the case of Weyl arrangements, it is proved by Orlik and Terao in \cite{OT2} that all restrictions are free, and such 
a free arrangement is called hereditarily free. Here let us check freeness of some restrictions 
of $\A$ by applying Corollary \ref{mille}.

Let 
$X_1=\{x_1=x_6\},\ 
X_2=\{x_1=x_6,x_2=x_5\},\ 
X_3=\{x_1=x_6, x_2=x_5,x_3=x_4\}$, and consider the freeness of 
$\A^{X_i}$ for $i=1,2,3$. 
Then it is easy to show that 
\begin{eqnarray*}
\A^{X_1}&:&x_6\prod_{2 \le i < j \le 6}(x_i^2-x_j^2))=0,\\
\A^{X_2}&:&x_5x_6\prod_{3 \le i < j \le 6}(x_i^2-x_j^2)=0,\\
\A^{X_3}&:&x_4x_5x_6\prod_{4 \le i < j \le 6}(x_i^2-x_j^2)=0.
\end{eqnarray*}
Since $\A^{X_3}$ is the Weyl arrangement of the type $B_3$, it is free with $\exp
(\A^{X_3})=(1,3,5)$. Hence we may apply Corollary \ref{mille} to check the freeness of these 
three arrangements.

By the equations, we can see that 
$|\A|=30,\ |\A^{X_1}|=21,\ 
|\A^{X_2}|=14,\ 
|\A^{X_3}|=9$. Hence Corollary \ref{mille} shows that $\A^{X_1}$ and $ \A^{X_2}$ are 
both free with $\exp(\A^{X_1})=(1,3,5,5,7)$ and
$\exp(\A^{X_2})=(1,3,5,5)$. 


\end{example}

\section{Supersolvable and divisionally free arrangements}

First recall the definition of the supersolvable arrangement.

\begin{define}
$\A$ is \textbf{supersolvable} if and 
only if there is a flag $\{X_i\}$ such that, 
$\A_{
X_i}$ is of rank $i$ for $i=0,\ldots,\ell-1$ and  
for every $H \neq L \in \A_{X_{i+1}} \setminus \A_{X_{i}}$, 
there is $K \in \A_{X_{i}}$ such that 
$H \cap L \subset K$. In this case, $\A$ is free with 
$\exp(\A)=(|\A|-|\A_{X_{\ell-1}}|,\ldots,|\A_{X_{2}}|-|\A_{X_{1}}|,
|\A_{X_{1}}|)$.
\label{SSdef}
\end{define}

Second, let us introduce a different definition of a supersolvable arrangement. 
We do not know whether it has been already known. Here we give a proof for the
completeness.

\begin{prop}
$\A$ is supersolvable if and only if there is a flag 
$\{X_i\}$ such that 
$$
b_2(\A)=\sum_{i=0}^{\ell-1} 
(|\A_{X_{i+1}}|-|\A_{X_{i}}|)|\A_{X_{i}}|.
$$
In this case, $\A$ is free with exponents 
$\exp(\A)=
(|\A_{X_\ell}|-|\A_{X_{\ell-1}}|,
|\A_{X_{\ell-1}}|-|\A_{X_{\ell-2}}|,\ldots,
|\A_{X_{2}}|-|\A_{X_{1}}|, 
|\A_{X_{1}}|)$.
\label{SS}
\end{prop}

\noindent
\textbf{Proof}. 
Let $\A_i:=\A_{X_i}$. Since a supersolvable arrangement is free with 
exponents in Definition \ref{SSdef}, the ``if'' part is immediate. Assume that 
$\A$ satisfies the equality in Proposition \ref{SS}. 
Assume that the assumption for supersolvable arrangements 
holds true for $\A_{0},\ldots,\A_{i}$. We show that, for 
any distinct $H,L \in \A_{i+1} \setminus \A_{i}$, there is $K \in \A_{i}$ such that 
$H \cap L \subset K$.

By the induction hypothesis, we know that $\A_{i}$ is supersolvable with 
$b_2(\A_{i})=\sum_{j=0}^{i-1}
(|\A_{j+1}|-|\A_{
j}|)|\A_{j}|$. Let $H \in \A_{i+1} \setminus \A_{i}$. Since 
$\A_{i} \cup \{H\}$ is of rank $i+1$ by definition of localization, 
it holds that 
$$
b_2(\A_{i} \cup \{H\}) \ge b_2(\A_{i})+|\A_{i}|.
$$
Hence 
\begin{eqnarray*}
b_2(\A_{i+1}) &\ge& b_2(\A_{i})+
(|\A_{i+1}|-|\A_{i}|)|\A_{i}|\\
&=&
\sum_{j=0}^{i}
(|\A_{j+1}|-|\A_{
j}|)|\A_{j}|.
\end{eqnarray*}
At every $i$, this has to be equal since we have the 
equation 
$$
b_2(\A) =
\sum_{j=0}^{\ell-1}
(|\A_{j+1}|-|\A_{
j}|)|\A_{j}|.
$$
Hence it holds that $
b_2(\A_{i} \cup \{H\}) = b_2(\A_{i})+|\A_{i}|\
$
for any $H \in \A_{i+1} \setminus \A_{i}$ and $i=0,\ldots,\ell-1$. 

Now assume that there are no $K \in \A_i$ such that $H \cap L \subset K$. 
Since $|(\A_i \cup\{H\})^H| \ge |\A_i
|$ by the definition of the localization, the above implies that $b_2(\A \cup \{H\}) >
b_2(\A_i) +|\A_i|$, which is a contradiction. Hence 
for any distinct $H, L \in \A_{i+1} \setminus \A_{i}$, there is $K \in \A_{i}$ such that 
$H \cap L \subset K$. \owari
\medskip

The reason why we introduced another characterization of 
supersolvable arrangements in Proposition \ref{SS}  is to point out the similarity of the supersolvablitiy to 
the divisional freeness introduced in \cite{A2}.

\begin{define}[Divisionally free arrangement, \cite{A2}, Definition 1.5]
$\A$ is divisionally free if  if there is a flag 
$\{X_i\}$ such that 
$$
b_2(\A)=\sum_{i=0}^{\ell-2} 
(|\A^{X_{i}}|-|\A^{X_{i+1}}|)|\A^{X_{i+1}}|.
$$
In this case, $\A$ is free with exponents 
$\exp(\A)=
(|\A^{X_0}|-|\A^{X_{1}}|,
|\A^{X_{1}}|-|\A^{X_{2}}|,\ldots,
|\A^{X_{\ell-2}}|-|\A^{X_{\ell-1}}|, 
|\A^{X_{\ell-1}}|)$. 
Such a flag is called a divisional flag.
\label{DF}
\end{define}

It is also shown in \cite{A2} that 
all inductively free arrangements are divisionally free (\cite{A2}, Theorem 1.6). Since 
supersolvable arrangements are inductively free, they are of course 
divisionally free. Here we give an another proof 
of the fact that supersolvable arrangements are divisionally free 
by using Proposition \ref{SS} and Definition \ref{DF} to see their similarity.

\begin{prop}[cf. \cite{A2}, Theorem 1.6]
A supersolvable arrangement $\A$ is divisionally free.
\label{SSDF}
\end{prop}

\noindent
\textbf{Proof}.
Let $\A$ be a supersolvable arrangement with a 
flag $\{X_i\}$ as in Definition \ref{SS}. Let 
$\alpha_1,\ldots,\alpha_\ell$ be linear forms such that 
$X_i=\{\alpha_1=\cdots=
\alpha_i=0\}$. Then 
define the flag
$\{Y_i\}$ by $Y_i:=\{\alpha_\ell=
\cdots=\alpha_{\ell-i+1}=0\}$. Then it is clear that this flag becomes a divisional 
flag.\owari

\begin{rem}
By Proposition \ref{SS}, Definition \ref{DF} and Proposition \ref{SSDF}, 
supersolvable and divisionally free arrangements are similar. They both use 
flags for localizations and restrictions respectively. Since there are a lot of nice properties 
for supersolvable arrangements, it is natural to ask whether 
some special properties which hold for supersolvable arrangements also hold true for 
divisionally free arrangements. 
\end{rem}

\end{document}